\documentclass[letterpaper, 10 pt, conference]{ieeeconf}  

\IEEEoverridecommandlockouts                              

\usepackage[latin1]{inputenc}
\usepackage{graphics} 
\usepackage{amssymb}  
\usepackage{tikz}
\usepackage{amsmath}

\title{\LARGE \bf 
  A discrete event traffic model explaining the traffic phases of the train dynamics on a linear metro
  line with demand-dependent control}

\author{Florian Schanzenbächer$^{1*}$, Nadir Farhi$^{2}$, Fabien Leurent$^{3}$, Gérard Gabriel$^{4}$
\thanks{$^{1}$ RATP, Paris, France and Université Paris-Est, France.}%
\thanks{$^{2}$ Université Paris-Est, COSYS, GRETTIA, IFSTTAR, France.}%
\thanks{$^{3}$ LVMT, ENPC, France. } %
\thanks{$^{4}$ RATP, Paris, France.}%
\thanks{* Corresponding author. {\texttt{florian.schanzenbacher@ratp.fr}}} }

\newtheorem{theorem}{Theorem}
\newtheorem{corollary}{Corollary}

\begin{document}

\maketitle
\thispagestyle{empty}
\pagestyle{empty}

\begin{abstract}
In this paper we present a mathematical model of the train dynamics in a linear metro line system
with demand-dependent run and dwell times.
On every segment of the line, we consider two main constraints.
The first constraint is on the travel time, which is the sum of run and dwell time.
The second one is on the safe separation time, modeling the signaling system,
so that only one train can occupy a segment at a time.
The dwell and the run times are modeled dynamically, with two control laws.
The one on the dwell time makes sure that all the passengers can debark from and embark into the train.
The one on the run time ensures 
train time-headway regularity in the case where perturbations do not exceed a run time margin.

We use a Max-plus algebra approach which allows to derive analytic formulas for the train time-headway and frequency depending
on the number of trains and on the passenger demand.
The analytic formulas, illustrated by 3D figures, permit to understand the phases of the train dynamics of a linear metro line being operated as a transport on demand system.
\end{abstract}


  \section{Introduction and Literature review}

Innovative public transportation systems are nowadays more and more turned towards transport on demand systems.
Transport on demand signifies that the transportation offer is adjusted in real-time to satisfy the passenger demand. 
This features trains stopping at stations just as long as needed to let passengers debark and embark.
Irregularities on the train time-headway
shall be recovered quickly to avoid a cascade effect, that is a longer headway causes a longer dwell time because of the accumulation of passengers, which causes an even longer headway. 
Therefore, transport on demand systems necessitate an efficient control. The latter detects and controls perturbations to offer a reliable service to costumers with a constant train frequency that equalizes the passenger charge over all trains. 



Several approaches for control of mass transit metro systems have been proposed.
The authors of~\cite{FCVC06} have applied a quadratic programming algorithm for optimal traffic control to a circular metro line, taking into
account a constant passenger arrival rate at the platforms.
In~\cite{SCF16} an optimal control strategy to gain on headway regularity for the traffic on a mass transit railway line using a quadratic programming algorithm has been proposed.
This method has been extended by~\cite{LSYG15} and~\cite{LDYG17}, considering uncertain time-variant passenger demand.
Moreover, the authors of~\cite{YIN16} have developed a real-time control approach for a metro system which deals with stochastic passenger demand on the one side,
and optimizes energy consumption, on the other side.
However, the models \cite{FCVC06},\cite{SCF16},\cite{LSYG15},\cite{LDYG17},\cite{YIN16} allow to represent only one train per inter-station~\cite{BCB91}, which limits their practical relevance, especially in networks
with longer inter-station distances.
The Max-plus approach features many advantages for the control of discrete event systems, like metro lines, 
which are discretized into segments. For more details on the Max-plus algebra refer to~\cite{BCOQ92}.

Demand-dependent control of transportation systems modeled in Max-plus algebra, has only recently begun to be treated in scientific
research.
The authours of~\cite{FNHL16a} have presented a Max-plus model for a linear metro line with constant dwell and run times. The authors of~\cite{FNHL16b} have adopted this model
for a linear metro line with dynamic dwell times depending on the passenger demand.
Their model includes an efficient control strategy which stabilizes the system in case of perturbations on the train time-headway.
In their approach, longer headways are recovered by shortening dwell times of retarded trains.
This control strategy is efficient since, in case of perturbations, trains running on schedule are not retarded, but trains running
behind schedule are accelerated via shorter dwell times to recover delays and to harmonize headways.
However, retarded trains generally face higher demand at the upcoming stops, because passengers have been accumulating over a longer time interval.
Consequently, recovering perturbations by cutting down dwell times is 	counterintuitive.

Therefore, in this paper we propose a mathematical Max-plus model for a linear metro line with demand-dependent dwell and run times.
The dwell time takes into account passenger upload and download rate as given by the vehicle characteristics, as well as passenger demand emissions and attractions of the platforms.
Perturbations on the train time-headway are recovered by shorter run times on the following inter-stations to avoid the amplification of the perturbation and a cascade effect.
The here presented mathematical model features a derivation of analytic formulas of the average train time-headway and frequency depending on the number of trains on the line and on the passenger demand.

In our future work, we will extend this model to cope with the stochastic nature of the passenger demand.
The authors of~\cite{LBFS12} have published on model predictive control for stochastic switching Max-plus linear systems.

\section{Review on Max-plus algebra}

As having been said above, Max-plus algebra has several advantages for modeling discrete event systems.

Max-plus algebra~\cite{BCOQ92} is the idempotent commutative semi-ring
%
($\mathbb{R} \cup \{ - \infty \}, \oplus, \otimes $),
where the addition and the multiplication of \textit{a, b} is defined by:
$a \oplus b = \max\{a,b\}$ and  $a \otimes b = a + b$.
This allows to model the dynamics of a metro
system in a matrix form which is highly advantageous for deriving analytic formulas of the
train dynamics, as the average train time-headway and
frequency, as well as for real-time traffic control.

For more details on how to model the train dynamics on a metro line in Max-plus algebra and how to technically
derive analytic formulas of the average train time-headway and the frequency 
refer to~\cite{FNHL16a}.
The main Theorem we use here to derive the analytic formulas for the average asymptotic train time-headway and frequency is from~\cite{CCGMQ98} and has already been applied in~\cite{FNHL16a}.

\section{Train dynamics}

\begin{figure}[htbp]
  \centering
  \includegraphics[scale=0.6]{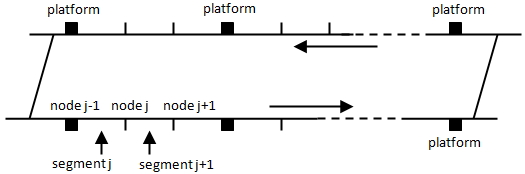}
  \caption{A linear metro line representation.}
  \label{fig-line}
\end{figure}

This work is based on~\cite{FNHL16a} where a discrete event traffic model explaining the traffic phases
of the train dynamics on a linear metro line has been proposed. In~\cite{FNHL16a}, lower bounds are imposed to
the run and dwell times, independent of the passenger travel demand. 
We extend here the model presented by the authors of~\cite{FNHL16a}, accounting for the passenger travel
demand in the run and dwell time dynamics.
As the authors of~\cite{FNHL16a}, we consider a space discretization of the metro line into $n$ segments,
as shown in Fig.~\ref{fig-line}.

Let us use similar notations to those used in~\cite{FNHL16a}.
\\~\\
\begin{tabular}{ll}
$d_{j}^{k}$ & the $k^{th}$ departure time of trains on segment $j$. \\
$a_{j}^{k}$ & the $k^{th}$ arrival time of trains on segment $j$.\\
  $r_{j}$ & the average running time of trains on segment $j$ \\
      & (between nodes $j-1$ and $j$). \\
  $w^k_{j}$ & $=d^k_{j}-a^k_{j}$ the $k^{\text{th}}$ dwell time on node $j$.
\end{tabular}  
\begin{tabular}{ll}
  $t^k_{j}$ & $=r_{j}+w^k_{j}$ the $k^{\text{th}}$ travel time from node \\
		& $j-1$ to node $j$.\\
  $g^k_{j}$ & $=a^k_{j}-d^{k-1}_{j}$ the $k^{\text{th}}$ safe separation time \\
                & (or close-in time) at node $j$.		
\end{tabular}  
\begin{tabular}{ll}		
  $h^k_{j}$ & $=d^k_{j}-d^{k-1}_{j} = g^k_{j}+w^{k}_{j}$ the $k^{\text{th}}$\\
      & departure time-headway at node $j$.\\
  $s^k_{j}$ & $=g^{k+b_{j}}_{j}-r_{j}$.\\    
  $b_{j}$ & $\in \{0,1\}$. It is $0$ (resp. $1$) if there is no \\
              & train (resp. one train) at segment $j$.\\
  $\bar{b}_{j}$ & $= 1 - b_{j}$.
\end{tabular}~~\\~~

Lower and upper bounds are written
$\underline{r}_j, \underline{w}_j, \underline{t}_j, \underline{g}_j, \underline{h}_j$
and $\bar{w}_j, \bar{t}_j, \bar{g}_j, \bar{h}_j$ where $\underline{w}_j$ is the minimum dwell time, regarding the demand, and $\underline{r}_j$ the minimum run time, given by the fastest speed profile. 
The average on $k$ and $j$ of the quantities above are denoted 
$r, w, t, g, h$ and $s$. We have the following relationships~\cite{FNHL16a}.
\begin{align}
    & g = r + s, \label{form1} \\
    & t = r + w, \label{form2} \\
    & h = g + w = r + w + s = t + s. \label{form3}
\end{align}    

We replace now the static run and dwell times by a model giving these two variables as a function of the demand and the vehicle characteristics.
Moreover, we suppose that the demand for the next train is completely satisfied by this train, that means all the passengers willing to debark will get off and
all the ones on the platform will board.
This assumption is necessary since, at this point, we do not take into account the capacity of platforms and trains.
It will be done in a further step.
Furthermore, the following passenger demand model supposes an average demand for each platform.
We use the following additional notations.
\\~\\
\begin{tabular}{ll}  
  $\lambda_{ij}$ & passenger travel demand from platform~$i$ to \\
      & platform~$j$, when $i$ and $j$ denote platforms, and\\
      & $\lambda_{ij} = 0$ if $i$ or $j$ is not a platform node.\\
  $\lambda^{\text{in}}_i$ & $=\sum_j \lambda_{ij}$ the average passenger arrival rate on \\
      & origin platform $i$ to any destination platform. \\      
  $\lambda^{\text{out}}_j$ & $=\sum_i \lambda_{ij}$ the average passenger arrival rate on \\
      & any origin platform to destination platform $j$. \\
  $\alpha^{\text{in}}_j$ & average passenger upload rate on platform $j$. \\      
  $\alpha^{\text{out}}_j$ & average passenger download rate on platform $j$. \\  
  $\frac{\lambda^{\text{out}}_j}{\alpha^{\text{out}}_j}\; h_j$ & $\approx \sum_{i} \lambda_{ij} h_{i} / \alpha_{j}^{out}$ time for passenger download\\
  & at platform $j$. \\
  $\frac{\lambda^{\text{in}}_j}{\alpha^{\text{in}}_j} \; h_j$ & time for passenger upload at platform $j$. \\
\end{tabular}  
\\~\\
We define $x_j$ as a passenger demand parameter
\begin{equation} \label{eq-x}
   x_j = \left( \frac{\lambda^{\text{out}}_j}{\alpha^{\text{out}}_j} + \frac{\lambda^{\text{in}}_j}{\alpha^{\text{in}}_j} \right),
\end{equation}
such that $x_j h_j$ represents the time needed for passenger down- and upload at platform $j$.

The model presented in~\cite{FNHL16a} considers two time constraints:
\begin{itemize}
  \item A constraint on the travel time on every segment $j$.
     \begin{equation}\label{const1}
       d^k_j \geq d^{k-b_{j}}_{j-1} + \underline{t}_j.
     \end{equation}
  \item A constraint on the safe separation time at every segment $j$.
     \begin{equation}\label{const2}
        d^k_j \geq d^{k-\bar{b}_{j+1}}_{j+1} + \underline{s}_{j+1}.
     \end{equation}
\end{itemize}

In~\cite{FNHL16a}, $\underline{t}_j$ in constraint~(\ref{const1}) is written $\underline{t}_j = \underline{r}_j + \underline{w}_j$
since running times $r_j$ are assumed to be constant over time, and dwell times $w_j$ are lower bounded without any
dependence on the passenger travel demand.
In the model we present here, we modify constraint~(\ref{const1}), whereas we keep constraint~(\ref{const2}) unchanged.
We replace constraint~(\ref{const1}) with the following.
\begin{equation}\label{const3}
   d^k_j \geq d^{k-b_{j}}_{j-1} + t^k_j(h^k_j,x_j) = d^{k-b_{j}}_{j-1} + r^k_j(h^k_j,x_j) + w^k_j(h^k_j,x_j),
\end{equation}
In constraint~(\ref{const3}), $r^k_j, w^k_j$ and thus $t^k_j$ are functions of the train time-headway $h^k_j$,
and of the passenger demand.
The function $t^k_j(h^k_j,x_j)$ is the control law which makes the train dynamics adaptive with respect to 
train delays and to passenger demand.

%

We pass here from static travel times, which are supposed to respect lower bounds, to a travel time model which is a function
of the passenger demand.

Moreover, the dwell time, the first component of the travel time, is a linear function of the train time-headway.
The run time, the second component of the travel time, is as well a function of the train time-headway and will be used
to recover perturbations on the dwell time.
We detail below the two functions $r^k_j(h^k_j)$ and $w^k_j(h^k_j)$.

\subsection{The dwell time control law}

Let us first recall that in~\cite{FNHL16b}, two dwell time models depending on
the passenger travel demand have been considered.
The first model is of the form
\begin{equation}\label{mod1}
  w^k_j = \frac{\lambda_j}{\alpha_j} h^k_j.
\end{equation}
We can see that the dwell time takes into account the needed upload time without taking into account the
needed download time. This model permits to adjust the dwell times on platforms in function of
the arrival demand.
However, it has been shown in~\cite{FNHL16b} that the train dynamics is not stable in this case.
That means train delays are amplified over time and propagate backwards through the metro line.

To deal with the latter problem, the authors of~\cite{FNHL16b} have proposed a second model for the train dwell times at
platforms. The model is of the form
\begin{equation}\label{mod2}
  w^k_j = \bar{w}_j - \theta_j \frac{\lambda_j}{\alpha_j} h^k_j,
\end{equation}
with $\theta_j$ being a control parameter to be fixed.
The authors have shown that the model~(\ref{mod2}) guarantees the stability of the train dynamics,
and that the dynamics admit an asymptotic regime with an asymptotic average train time-headway.
The latter is derived by simulation in function of the number of trains and of the level of
the passenger demand.
However, with the control law of  model~(\ref{mod2}), a delayed train at platform $j$, which induces an accumulation of passengers at platform $j$, will reduce the train dwell time at that platform,
what is absurd for passengers.

We propose here a model which resolves all those problems. 
The dwell time model is similar to the model~(\ref{mod1}), but it also takes into account
the attraction term of the travel demand.
In order to deal with the instability issue, we complete the dwell time model with a running time
model which cancels the terms that causes instability in the dwell time model.
More precisely, for a train delay at a platform $j$, the dwell time model extends the dwell time
at that platform in order to satisfy the needed download and upload times.
The running time model will reduce the running time from platform $j$ to platform $j+1$
in order to compensate in such a way that the whole travel time ($t=w+r$) is stable.

The dwell time model is the following.
\begin{align}
  & w_j^k(h_j^k,x_j) = \min (x_j h_j^k, \bar{w}_j ), \label{eq-dwell}
\end{align}

where the $k^{th}$ dwell time on platform $j$ is the minimum time for passenger up- and download
on the platform, and a maximum dwell time $\bar{w}_j$ to avoid train congestion behind a train stopping too long. 
The lower bounds on $w_j$ are fixed as follows (upper bounds are fixed accordingly).
\begin{align}
  & \underline{w}_j = x_j \underline{h}_j,
\end{align}
where $\underline{h}_j$ are derived from given $\underline{g}_j$, from the formula 
$h = g + w$ of~(\ref{form3}), and from the dwell time law~(\ref{eq-dwell}), as follows. We have
\begin{align}
   & \underline{h}_j = \underline{g}_j + \underline{w}_j = \underline{g}_j + x_j \underline{h}_j.
\end{align}
Then we get
\begin{align}
   & \underline{h}_j = 1/(1-x_j) \; \underline{g}_j,  
\end{align}
and then
\begin{align}
   & \underline{w}_j = X_j \underline{g}_j, \label{eq-min-dwell}
\end{align}
where $X_j := x_j / (1 - x_j)$.
Let us notice, that~(\ref{eq-min-dwell}) corresponds to the minimum dwell time satisfying the demand.
The upper bounds are fixed similarly.
From the point on we observe a small perturbation on the train time-headway, the following longer
dwell time, modeled in equation~(\ref{eq-dwell}), will make this perturbation growing,
since longer headways result in longer dwell times which results in even longer headways.
Our aim is to recover perturbations on the dwell time by a dynamic running time model in order to stabilize the system.

\subsection{The running time control law}

We propose the following running time law.
\begin{align}
  & r_j^k(h_j^k,x_j) = \max \left\{ \underline{r}_j, \tilde{r}_j - x_j \left( h_j^k - \underline{h}_j \right) \right\}, \label{eq-run}
\end{align}
where $\tilde{r}_j$ is the average (nominal) running time of trains on segment $j$.

The model~(\ref{eq-run}) gives the running time as the maximum between a given minimum running time $\underline{r}_j$
and a term that subtracts $x_j \left( h_j^k - \underline{h}_j \right)$ from the nominal running time.
The term $x_j \left( h_j^k - \underline{h}_j \right)$ expresses a deviation of the upload and download time, 
due to a deviation of the train time-headway.
We notice here that the term $x_j h^k_j$, appearing in the dwell time law~(\ref{eq-dwell}) with a sign ``$+$'',
appears in the running time law~(\ref{eq-run}) with a sign ``$-$''.



Combining the dwell time law~(\ref{eq-dwell}) with the running time law~(\ref{eq-run}), we obtain the
following train travel time law.
\begin{equation}\label{eq-travel}
   t_j^k(x_j) = r_j^k(h_j^k,x_j) + w_j^k(h_j^k,x_j).
\end{equation}

Let us use the notations.
\begin{align}
    & \Delta h_j := \bar{h}_j - \underline{h}_j, \quad \Delta g_j := \bar{g}_j - \underline{g}_j, \nonumber \\
    & \Delta w_j := \bar{w}_j - \underline{w}_j, \quad \Delta r_j := \tilde{r}_j - \underline{r}_j, \nonumber
\end{align}
with $\bar{h}_j$ being the longest headway observed.   
It is then easy to check the following.
$$\Delta w_j = x_j \Delta h_j = X_j \Delta g_j, \forall j.$$
Then we have the following result.
\begin{theorem}\label{thm-mp}
  If $h^1_j \leq \bar{h}_j=1/(1 - x_j) \; \bar{g}_j, \forall j$ and $\Delta r_j \geq \Delta w_j = X_j \Delta g_j, \forall j$, then the dynamic system~(\ref{const2})-(\ref{const3}) is
  a Max-plus linear system, and is equivalent to
  \begin{align}
     & d^k_j \geq d^{k-b_{j}}_{j-1} + \tilde{r}_j + X_j \underline{g}_j. \label{eq-mp1} \\
     & d^k_j \geq d^{k-\bar{b}_{j+1}}_{j+1} + \underline{s}_{j+1}. \label{eq-mp2}
  \end{align}
\end{theorem}
~~\\~~
\proof 
By induction, let us first show that~(\ref{const2})-(\ref{const3}) is equivalent to~(\ref{eq-mp1})-(\ref{eq-mp2}) for $k=1$.
\begin{itemize}
  \item On the one side, we have $x_j h^1_j \leq x_j \bar{h}_j = \bar{w}_j$. Therefore, the first term realizes the minimum in~(\ref{eq-dwell}).
    That is,~(\ref{eq-dwell}) is equivalent to
    \begin{equation}\label{eq-dwell1}
       w_j^1(h_j^1,x_j) = x_j h_j^1.
    \end{equation}
  \item On the other side, we have
    $$\begin{array}{ll}
            \tilde{r}_j - x_j(h^1_j - \underline{h}_j)  & \geq \tilde{r}_j - x_j\left( \bar{h}_j - \underline{h}_j \right)\\
                                                        & = \tilde{r}_j - \Delta w_j \geq \tilde{r}_j - \Delta r_j = \underline{r}_j.
      \end{array}$$
    Therefore, the second term realizes the maximum in~(\ref{eq-run}).
    That is,~(\ref{eq-run}) is equivalent to
    \begin{equation}\label{eq-run1}
       r_j^1(h_j^1,x_j) = \tilde{r}_j - x_j \left( h_j^1 - \underline{h}_j \right).
    \end{equation}
\end{itemize}
Consequently,~(\ref{eq-travel}) gives
\begin{equation}\label{eq-travel1}
   t_j^1(x_j) = \tilde{r}_j + x_j \underline{h}_j = \tilde{r}_j + X_j \underline{g}_j.
\end{equation}
Then~(\ref{const3}) can be written
\begin{equation}\nonumber
   d^1_j \geq d^{1-b_{j}}_{j-1} + \tilde{r}_j + X_j \underline{g}_j.
\end{equation}

Let us now show, that if~(\ref{const2})-(\ref{const3}) is equivalent to~(\ref{eq-mp1})-(\ref{eq-mp2}) for a given $k$, then
it holds also for $k+1$.
Since it holds for $k$, then a Max-plus linear dynamics~(\ref{eq-mp1})-(\ref{eq-mp2}) will be applied for $k$.
We notice by $\mathbf f$ the Max-plus map of the Max-plus dynamics.
We know that Max-plus linear maps are $1$-Lipschitz for the sup. norm.
The assertion holds for $k$ means that $h^k_j \leq \bar{h}_j, \forall j$.
Then $||d^k - d^{k-1}||_{\infty} \leq \bar{h}_j$.
Hence
$$||d^{k+1} - d^{k}||_{\infty} = ||\mathbf f(d^k) - \mathbf f(d^{k-1})||_{\infty} \leq ||d^k - d^{k-1}||_{\infty} \leq \bar{h}_j.$$
Therefore $h^{k+1}_j \leq \bar{h}_j, \forall j$.
Then, we can easily show (as done for $k=1$) that~(\ref{const2})-(\ref{const3}) is equivalent to~(\ref{eq-mp1})-(\ref{eq-mp2}) for $k+1$.
\endproof

Let us interpret the two conditions of Theorem~\ref{thm-mp}.
\begin{itemize}
 \item Condition $h^1_j \leq \bar{h}_j=1/(1 - x_j) \; \bar{g}_j, \forall j,k$ limits the initial headway $h^1_j$
   (i.e. the initial condition) to 
   its upper bound $\bar{h}_j$, which is given by the level of the passenger travel demand $x_j$ at platform $j$,
   and by the upper bound $\bar{g}_j$ on the safe separation time at the same platform.
   We notice here that a big value of $h^k_j$ corresponds to a delay of the $k^{th}$ train passing by platform $j$.
   Therefore, this condition, tells that if all the delays expressed by $h^k_j, \forall j$ are limited to $\bar{h}_j, \forall j$,
   then the dynamic system is Max-plus linear, and is then stable, and admits a stationary regime.
   In other words, the train dynamics is stable under small disturbances.
 \item Condition $\Delta r_j \geq \Delta w_j = X_j \Delta g_j, \forall j$ limits the margin on the train dwell times
   to the margin on the train running times. 
   Indeed, as explained above, the model consists in responding to disturbances and train delays by first extending
   the train dwell times so that the passengers accumulated on the platforms have time to embark, and second,
   by recovering the dwell time extension by reducing the train run times on the downstream inter-stations.
   In order to ensure this recovering to be possible, and then the dynamic system to be stable, we need the margin
   on the running times to be at least equal to the margin on the train dwell times at the platforms.
\end{itemize}

It is important to notice that fixing bounds on $r$ and on $g$ will impose bounds on $s$, according to~(\ref{form1}).
$$\underline{g}_j = \underline{r}_j + \underline{s}_j.$$


In~(\ref{eq-mp1}), the travel time $t^k_j(x_j)$ depends on the demand $X_j$,
weighted by the minimum safe separation time $\underline{g}_j$.
This dependence is linear on $X_j$.
We notice also that, in general we have
$$\lambda^{\text{out}}_j, \lambda^{\text{in}}_j << \alpha^{\text{out}}_j, \alpha^{\text{in}}_j.$$
Therefore $x_j$ can be assumed to be close to $0$.
By consequent, a linearity with respect to $X_j$ is a kind of exponential relationship with respect to $x_j$,
since $X_j = x_j/(1 - x_j)$.
%
%
%
%
Condition $\Delta r_j \geq \Delta w_j = X_j \Delta g_j, \forall j$ of Theorem~\ref{thm-mp} makes a link between two important parameters
which are the margin on the running time control $\Delta r_j$,
and the level of the travel demand $X_j$.
%
%
%
%
Finally, the dynamic modeling of the dwell and run time allows to adjust the run time margin $\Delta r_j$ depending on the demand.
The system can be optimized with regard to train frequency or stability.

In the following section, we give directly the analytic formulas for the average asymptotic train time-headway $h$ and the average asymptotic
frequency $f$ for the three traffic phases, which can be derived as in~\cite{FNHL16a}.

\section{The traffic phases}

The Max-plus theorem derived in~(\cite{FNHL16a}, Theorem 1) gives the average asymptotic train time-headway on a linear metro line as a function of the number of trains circulating on the line.
Let us consider the following notations.
\\~\\
\begin{tabular}{ll}  
  $m$ & the number of trains.\\
  $n$ & the number of segments on the line.\\  
  $X$ & the vector of demand parameters $X_j, \forall j$.
\end{tabular}  
\\~\\
Since the here presented model can, as the model with static run and dwell times in~\cite{FNHL16a}, be written in Max-plus algebra,
we can directly replace $\underline{t}_j$ in~(\cite{FNHL16a}, Theorem 2) by $t_j(x_j)$.

\begin{equation}\label{th-h}
  h (m) = \max \left\{ \frac{\sum_j t_j(x_j)}{m}, \max_j (t_j(x_j)+\underline{s}_j), \frac{\sum_j \underline{s}_j}{n-m} \right\}.
\end{equation}

By replacing $t^k_j(x_j)$ using~(\ref{eq-travel}), we obtain the following result.

\begin{theorem}\label{thm-1}
The average asymptotic train time-headway of the linear Max-plus system with dynamic dwell and run times depending on the demand is given by the average asymptotic growth rate
of the system.
The average headway depends on the number of trains $m$ on the line and the passenger travel demand for every station $X_j$.
   $$  h(m,X) = \max \left\{ \begin{array}{l}
                               \frac{\sum_j (\underline{g}_j X_j + \tilde{r}_j)}{m}, \\~\\
                               \max_j ((\underline{g}_j X_j + \tilde{r}_j)+\underline{s}_j),\\~\\
                               \frac{\sum_j \underline{s}_j}{n-m}.
                            \end{array} \right.$$                                            
\end{theorem}
~~\\~~
\begin{proof}
Under the conditions of Theorem~\ref{thm-mp}, the system can be written in Max-plus algebra.
It has been shown in~\cite{FNHL16a} that in this case, the average asymptotic growth rate of the system can be analytically derived and
corresponds to the average headway.
We replace $\underline{t}_j$ in the Max-plus theorem of~\cite{FNHL16a} and obtain demand-dependent headway equation.
\end{proof}
We notice that $h$ depends not only on the number of trains $m$, but furthermore on a weighted mean of the demand
$\sum_j (\underline{g}_j X_j + \tilde{r}_j)$, and on the maximum of all the station demands $\max_j ((\underline{g}_j X_j + \tilde{r}_j)+\underline{s}_j)$.

We obtain the analytic formulas for the average train frequency on the line directly from Theorem~\ref{thm-1}.
\begin{corollary}\label{cor-1}
The average frequency of the linear Max-plus system with demand-dependent dwell and run times is a function of the number of trains and the travel demand.
	$$f (m,X) = \min \left\{ \begin{array}{l}
					\frac{m}{\sum_j (\underline{g}_j X_j + \tilde{r}_j)}, \\~\\
					\frac{1}{\max_j ((\underline{g}_j X_j + \tilde{r}_j)+\underline{s}_j)}, \\~\\
					\frac{n-m}{\sum_j \underline{s}_j}.
					 \end{array} \right.$$
\end{corollary}
~~\\~~
\begin{proof}
Directly from Theorem~\ref{thm-1} with $f = 1/h$.
\end{proof}

\begin{table*}[thbp]
\centering
\vspace{10pt}
\caption{Average asymptotic train time-headway and frequency in function of number of trains and demand (x/ X). Three traffic phases to be distinguished: headway
1-linear in \textit{X} for \textit{m = const}; 2-linear in \textit{X}; 3-independent of \textit{X}. Linear metro line, Data RATP, Paris.}
\begin{tabular}{|c|c|}
  \hline
  $h(m,x)$ & $h(m,X)$ \\
  \hline
   & \\
  \includegraphics[scale=0.20]{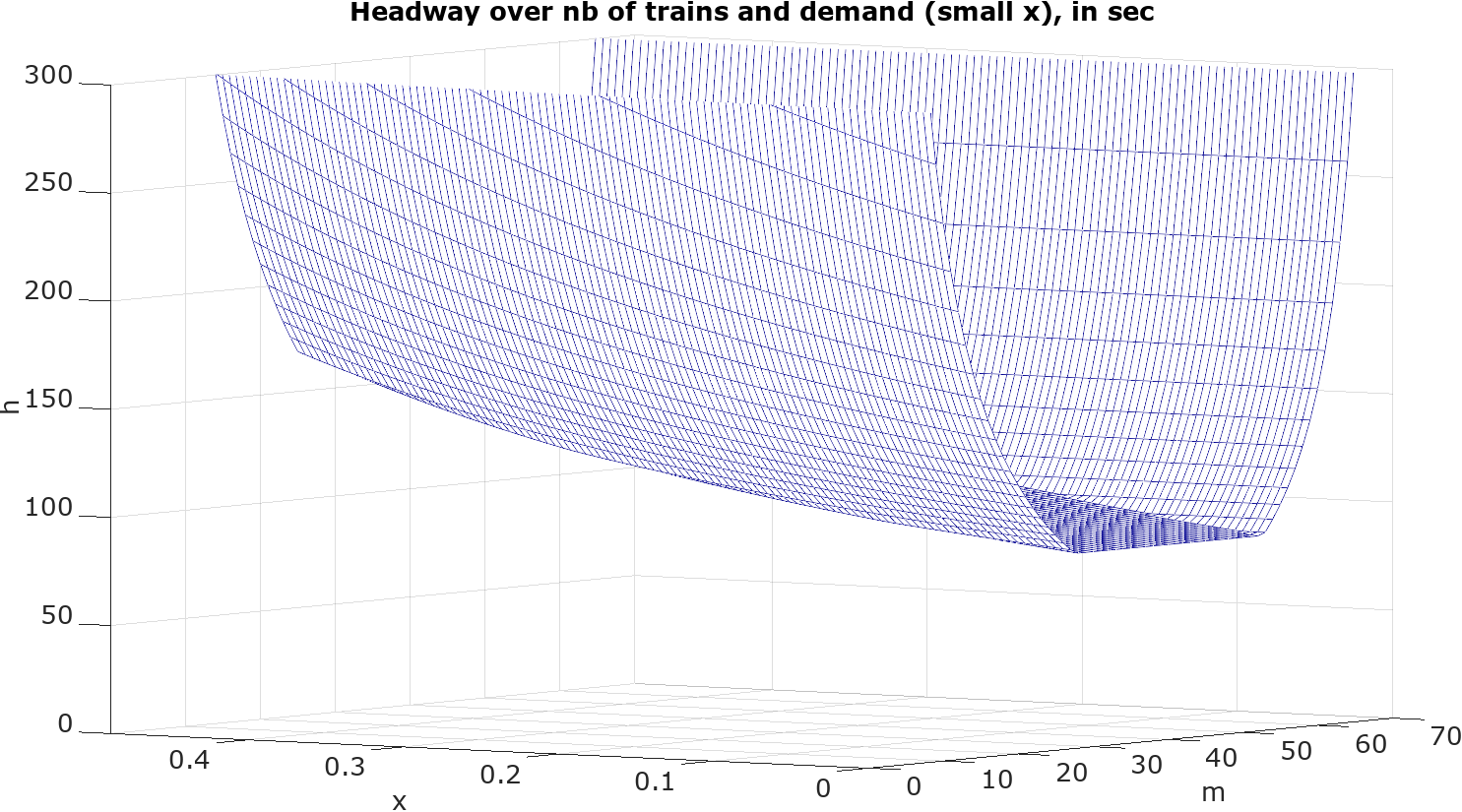} & \includegraphics[scale=0.19]{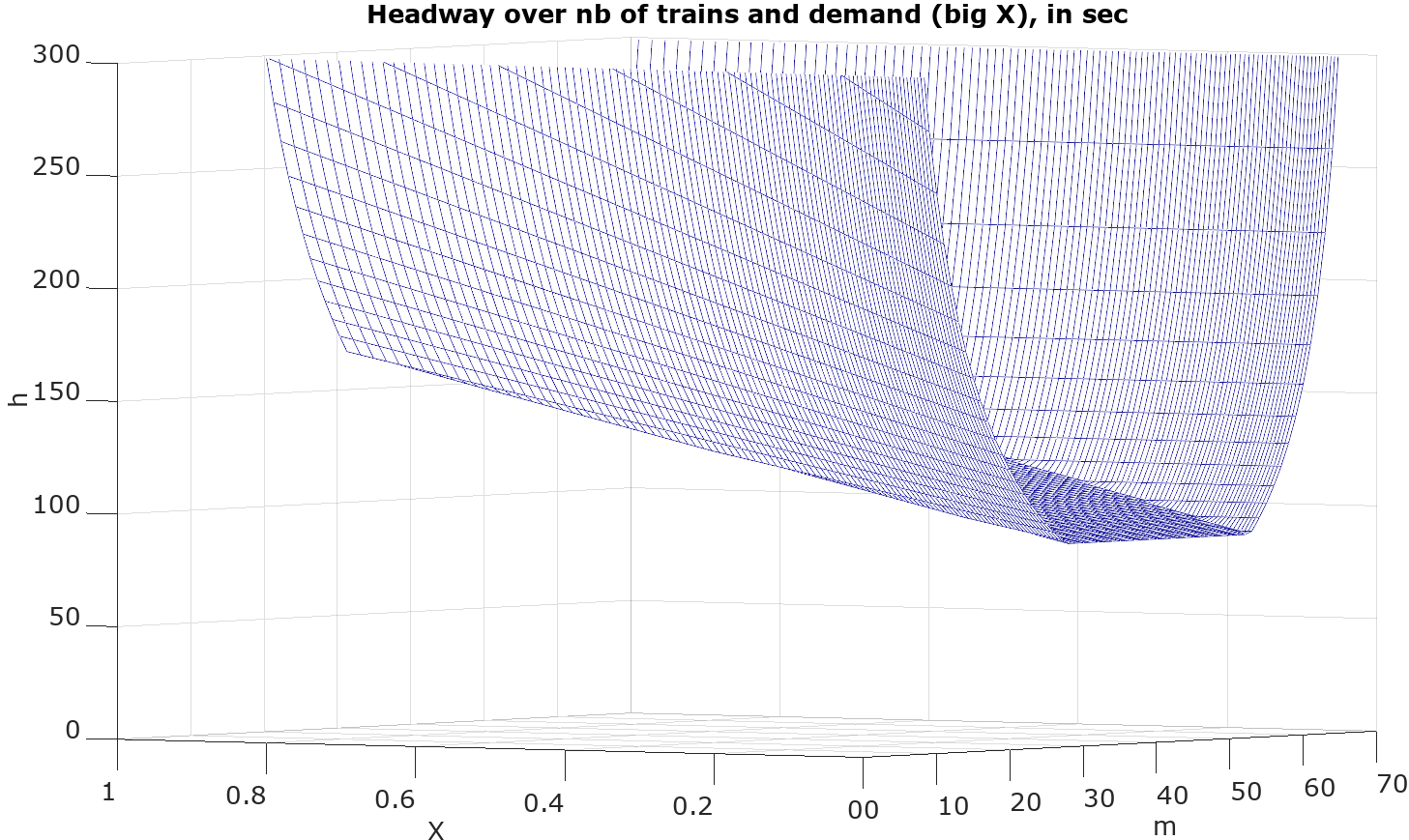} \\
  \hline 
  $f(m,x)$ & $f(m,X)$ \\
   \hline
   & \\
  \includegraphics[scale=0.20]{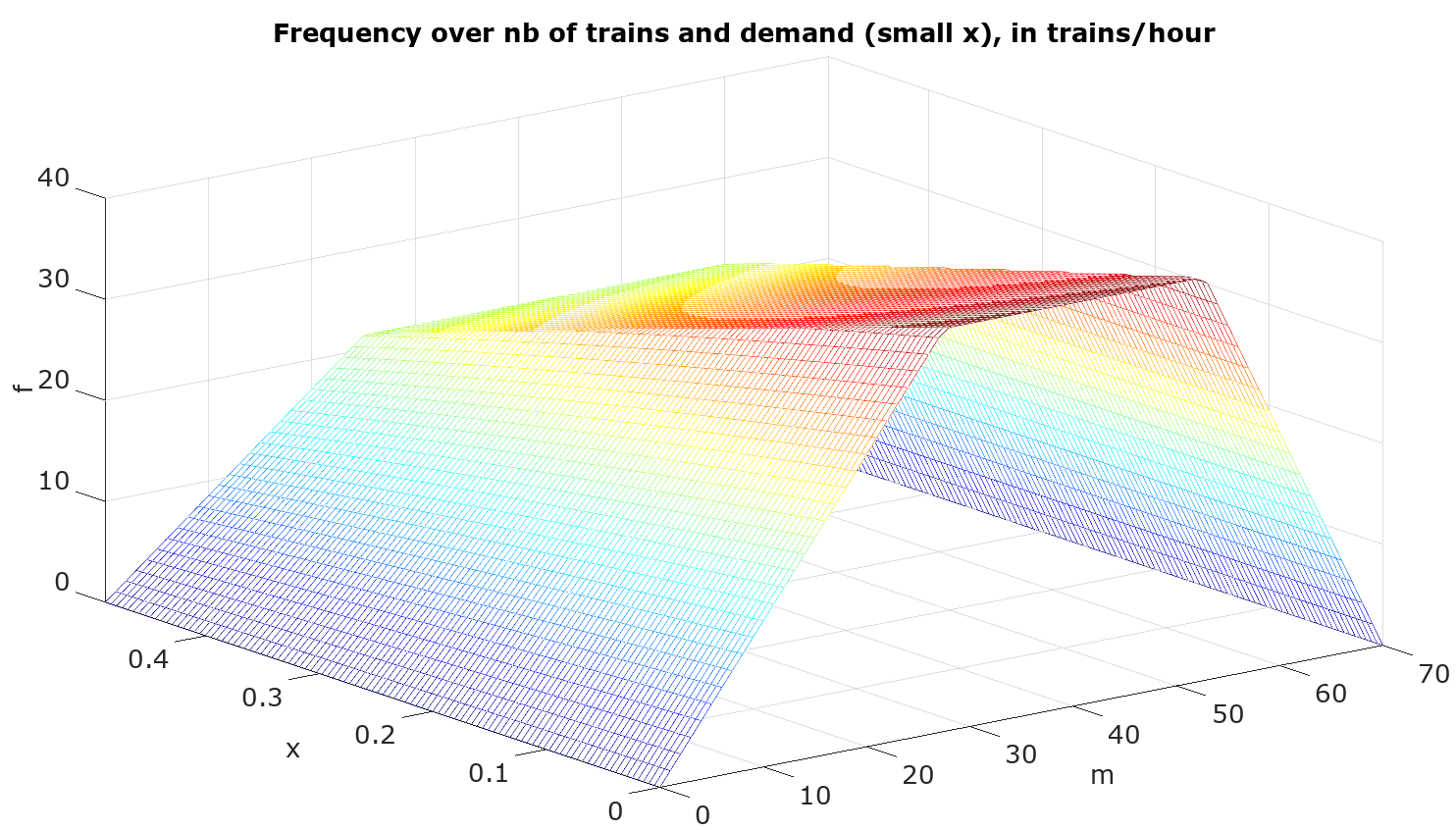} & \includegraphics[scale=0.20]{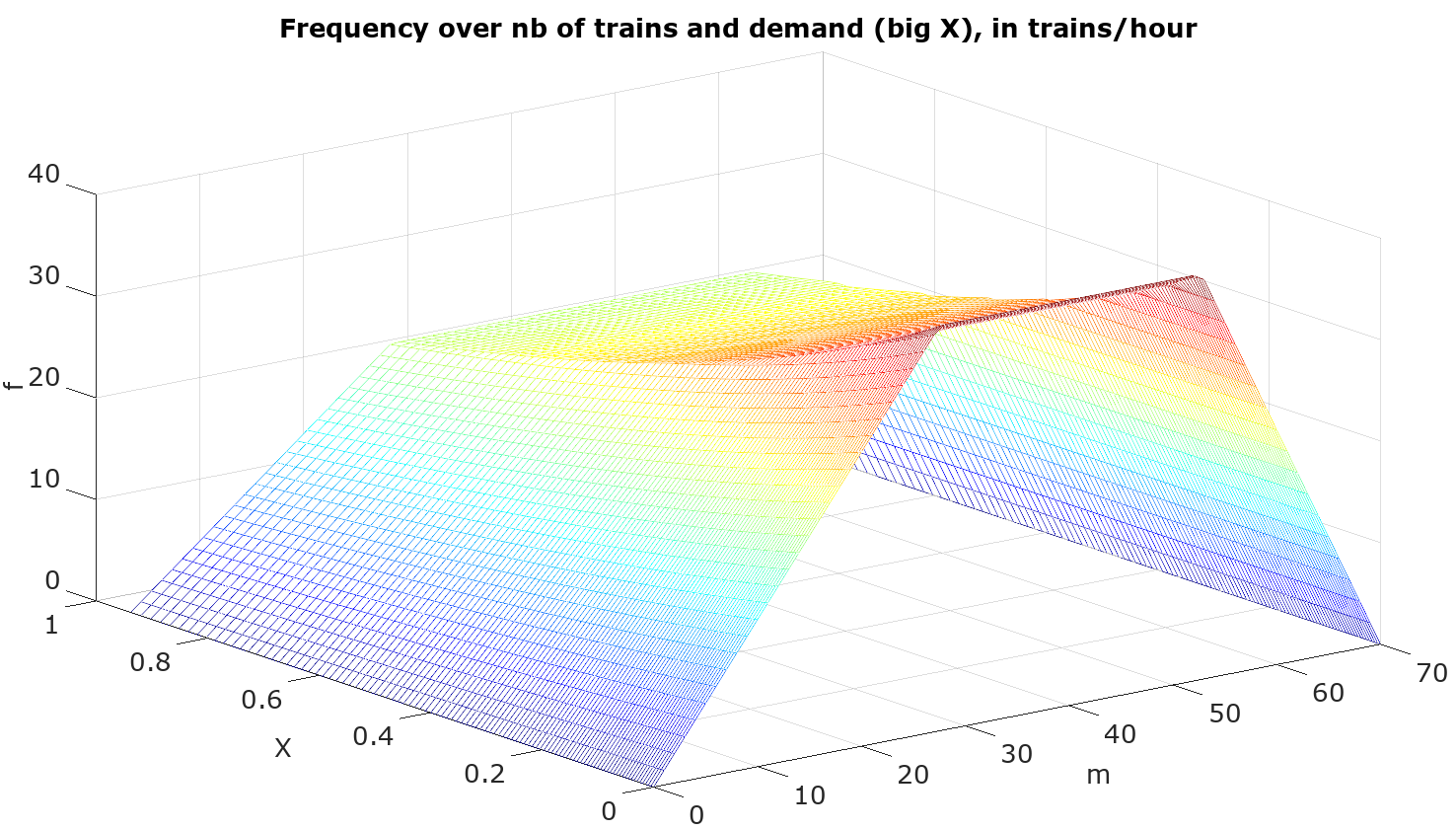} \\
  \hline 
  & \\
  \includegraphics[scale=0.20]{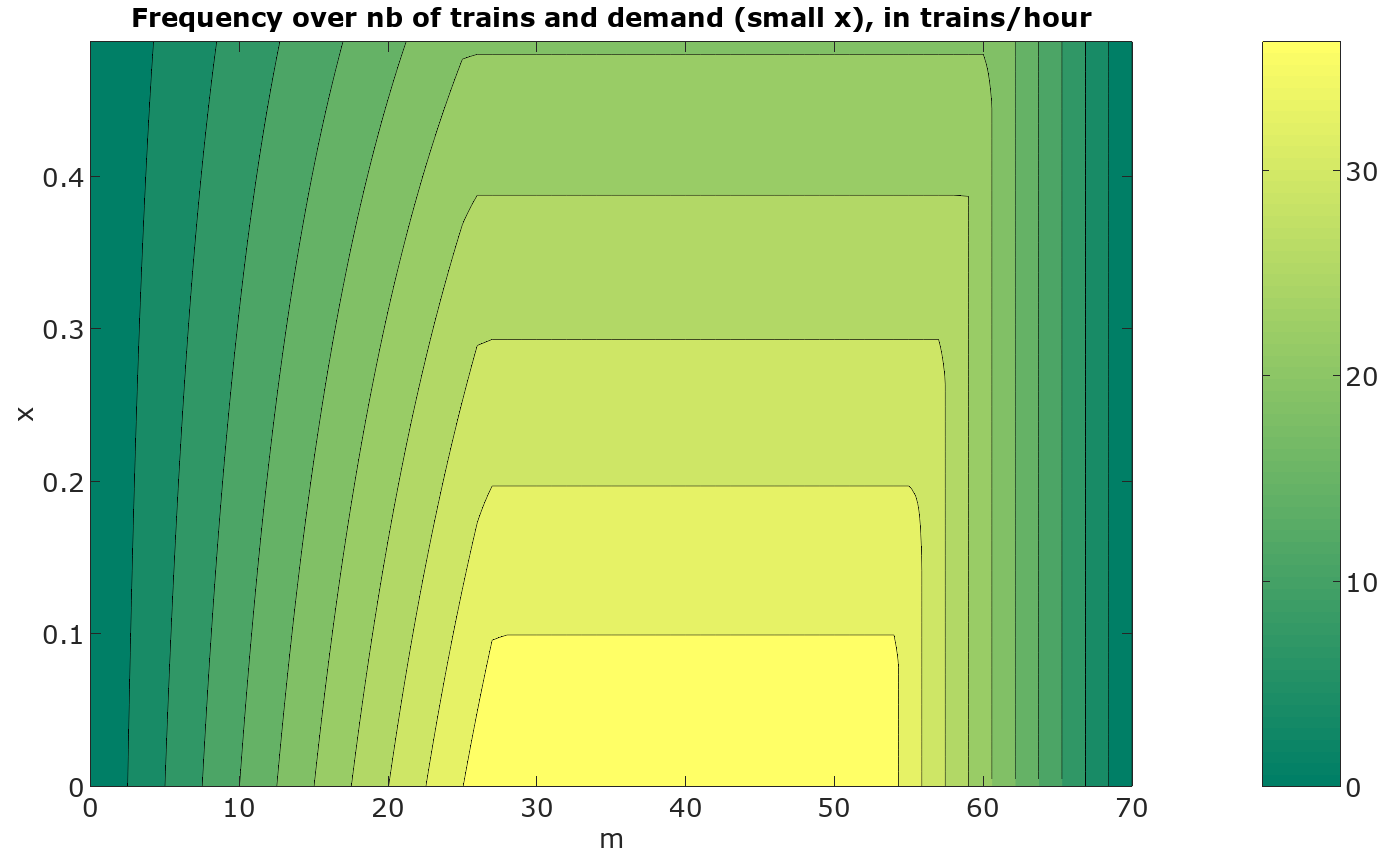} & \includegraphics[scale=0.20]{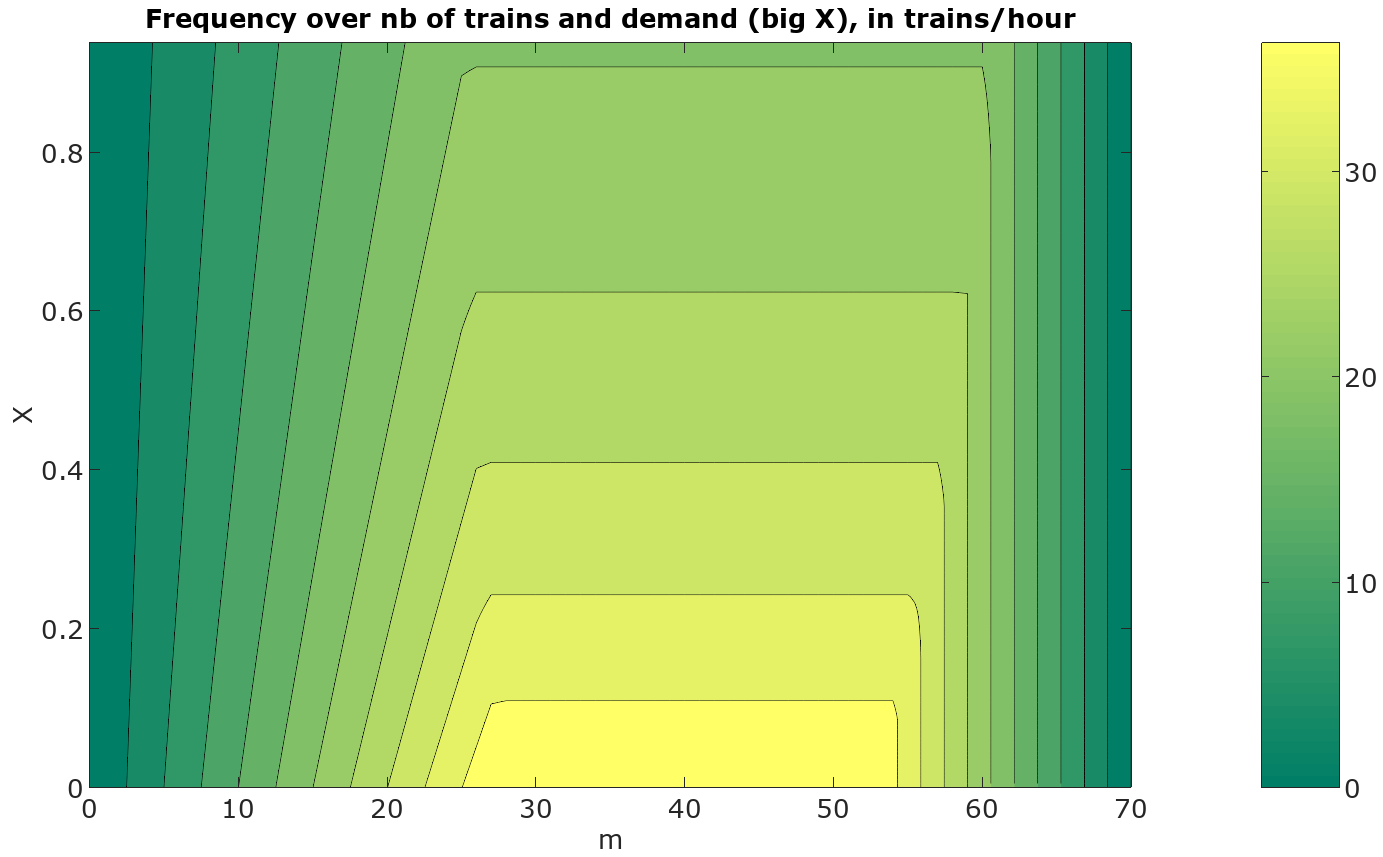} \\
  \hline
\end{tabular}
\label{tab}
\end{table*}

Table~\ref{tab} depicts the three headway phases derived in Theorem~\ref{thm-1}, as well as the three frequency phases of Corollary~\ref{cor-1} for a linear metro line with 18 stations from the RATP network in Paris. The demand level has been assumed to be equal for all the platforms.
The graphs in the left column are given in function of the average demand parameter $x = (\lambda^{out}/\alpha^{out} + \lambda^{in}/\alpha^{in})$ 
, whereas the ones in the right column depict the main results over $X = x/(1-x)$.

Let us interpret the first two graphs with the results of Theorem~\ref{thm-1}.
The headway in the free flow traffic phase depends on $X$ and on $m$. For a fixed $m$, the headway is linear in $X$.
In the maximum frequency phase, the headway is given by the maximum $\max_j{((\underline{g}_jX_j + \tilde{r}_j) + \underline{s}_j)}$.
The headway in the maximum frequency phase is here linear in $X$ and moreover independent of $m$.
In the congested traffic phase, the analytic formula shows that the headway is independent of the demand.
Furthermore, the headway increases with $m$.

On the left side, the headway is depicted over $x$.
In the free flow traffic phase, the headway is non linear, even for a fixed $m$.
Since $X = x/(1-x)$, the headway is non linear in the maximum frequency phase.
Finally, in the congested traffic phase, the headway is independent of the demand and increases with $m$.

Let us now look at the frequency, as given by Corollary~\ref{cor-1}. 
In the free flow traffic phase, for a fixed $m$, the frequency is non linear in $X$, since the headway is linear here and $f = 1/h$.
For a fixed $X$, the frequency is linear in $m$.
The same can be observed for the maximum frequency phase, which is non linear in $X$, because the headway is linear in $X$ here.
In the congested traffic phase, the frequency is linear in $m$.

For the frequency given over $x$, the free flow traffic phase in non linear for fixed $m$, but it is linear in $m$ for a fixed $x$.
In the maximum frequency phase, the frequency is independent of the number of trains $m$, and it is non linear in the demand parameter $x$.
Finally, as above, the frequency is linear in $m$ in the congested traffic phase, whereas it degrades linearly.
%

\section{CONCLUSIONS}

We have presented here a mathematical model which describes the train dynamics of a linear metro line, and
proposes a traffic control strategy.
The model considers two control laws, one for the train dwell and one for the run times, taking into account the
passenger demand.
It is a discrete event traffic model.
The approach is based on the Max-plus algebra which has allowed us to derive analytic formulas
for the average asymptotic
train time-headway and frequency, in function of the number of trains and of the passenger demand,
with the characteristics of the infrastructure and the vehicles being taken into account, too.

Our approach solves the problem 
that delayed trains facing a higher demand on the upcoming stations, are more and more delayed, 
which leads to a serious destabilization of the system. 
The model detects such a possible cascade effect.
In case of a delayed train to a given platform, where passengers have been accumulating,
the train dwell time is extended to respond to the passenger demand.
To avoid the amplification of the train delay, the running time is reduced in the following
inter-station to recover the dwell time extension.
Consequently, the train dynamics responds to the passenger demand, and remains stable.

Theorem~\ref{thm-mp} has shown that dwell times can only be extended up to the margin on the running time
in the following inter station. However, as given in Theorem~\ref{thm-1}, the average asymptotic headway
increases with the margin on the run time. That means, a compromise has to be found to satisfy both a
high frequency service on the one side, and to guarantee system stability, on the other side.

In our future work, we will apply this approach to a metro line with a junction. Moreover, we will propose a model for the passenger stock at the platforms and in the trains. In an ultimate version, we
will finally handle the uncertainties in the system, due to the stochastic nature of the passenger demand.

%
%
%
%
%
%
%
%
%
%




\end{document}